\documentclass[12pt]{article}
\usepackage[cp1251]{inputenc}
\usepackage[english,russian]{babel}
\usepackage{amsmath}
\usepackage{amssymb}
\usepackage{amsfonts}
\usepackage{graphicx}

\newtheorem{lemma}{Lemma}
\newtheorem{theorem}{Theorem}

\newtheorem{proposition}{Proposition}

\begin{document}

\renewcommand{\refname}{References}
\renewcommand{\figurename}{Fig.}
\renewcommand{\tablename}{Table}
\renewcommand{\proofname}{Proof.}
\thispagestyle{empty}

\begin{center}
\textbf{\Large{On Jacobian group and complexity of the $Y$-graph}}

\vspace{0.5cm}
\textbf{Y.S. Kwon, A.D. Mednykh, I.A. Mednykh}%

\end{center}

{\small\begin{quote}
\noindent{\sc Abstract.} In the present paper we suggest a simple approach for counting Jacobian group of the $Y$-graph $Y(n; k, l, m).$ In the case $Y(n; 1, 1, 1)$ the structure of the Jacobian group will be find explicitly. Also, we obtain a closed formula for the number of spanning trees of $Y$-graph in terms of Chebyshev polynomials and give its asymtotics.
\medskip

 \end{quote}}

\sloppy

\section{Introduction}

The notion of the Jacobian group of a graph, which is also known as the Picard group, the critical group, and the dollar or sandpile group, was independently introduced by many authors (\cite{CoriRoss}, \cite{BakerNorine}, \cite{Biggs}, \cite{BachHarpNagnib}).
This notion arises as a discrete version of the Jacobian in the classical theory of Riemann surfaces. It also admits a natural interpretation in various areas of physics, coding theory, and financial mathematics.
The Jacobian group is an important algebraic invariant of a finite graph. In particular, its order coincides with the number of spanning trees of the graph,
which is known for some simple graphs such as the wheel, fan, prism, ladder, and M\"oebius ladder \cite{BoePro}, grids \cite{NP04}, lattices \cite{SW00}, Sierpinski gaskets \cite{CCY07,AD11}, 3-prism and 3-anti-prism \cite{SWZ16}  and so on.
 At the same time, the structure of the Jacobian is known only in particular cases \cite{CoriRoss}, \cite{Biggs}, \cite{Lor}, \cite{YaoChinPen}, \cite{ChenHou}, \cite{MedZind}
and \cite{MedMed2}.

We mention that the number of spanning trees for circulant graphs is expressed in terms of the Chebyshev polynomials; it was found in \cite{ZhangYongGol}, \cite{ZhangYongGolin}, \cite{XiebinLinZhang}, \cite{MedMed17}  and \cite{MedMed18}.
More generally, this result also holds for arbitrary cyclic coverings of a graph \cite{KwonMedMed20}. In particular, (see \cite{KMM} and \cite{Ilya}) this is true for the generalized Petersen graph $GP(n, k)$  and $I$-graph  $I(n,j,k).$
These two graphs are expansions of the tree consisting of a single edge. In the same time, there is a wide family of cubic graphs that are expansions of the tree consisting of more than one edge. Three families of such graphs were discovered by Biggs \cite{BigIYH}.
He calls these graphs the $I$-graph, $Y$-graph, and $H$-graph, because of the method by which the graphs were created. This method of construction was generalized in \cite{HorBou}.
In the present paper we investigate the structure of Jacobian and find the number of spanning trees for a $Y$-graph. The similar questions for $H$-graphs will be considered in the forthcoming paper by the authors.

A $Y$-graph $Y(n; k, l, m)$ has $4n$ vertices arranged in four segments of $n$ vertices. Let the vertices be $v_{x,y}$ for $x = 0, 1, 2, 3,$ and $y $ in the integers modulo $n.$
Call the graph induced by the $v_{x,y}$ for a fixed $x,$ the $x$-{\it segment.} The $0$-segment is called the $\it inner \, segment$ and consists of $n$ independent vertices; while each of the other three segments is called an {\it outer segment.}
These three segments consist of {\it outer} edges which join each $v_{1,y}$ to $v_{1,y+k},$ each $v_{2,y}$ to $v_{2,y+l},$ and each $v_{3,y}$ to $v_{3,y+m},$ where the subscript addition is performed modulo $n.$
Each outer vertex $v_{i,y}$ is joined to the inner vertex $v_{0,y}$ by an {\it inner} edge. Thus $Y(n; k, l, m)$ consists of $n$ inner vertices, $3n$ outer vertices, $3n$ inner edges, and $3n$ outer edges.
For this definition to produce a connected graph it is necessary that $\gcd(k, l, m, n)=1$, where $\gcd(k, l, m, n)$ is the greatest cmmon divisor of $k, l, m$ and $n$. In what follows, we restrict ourselves by connected graphs only.

In this paper, we find an approach for counting Jacobian group of the $Y$-graph $Y(n;k,l,m).$ It the case $Y(n; 1,1,1)$ the structure of the Jacobian group will be find explicitly. Also, we obtain a closed formula for the number of spanning trees of $Y$-graph in terms of Chebyshev polynomials.

\section{Basic definitions and preliminary facts}

Consider a connected finite graph $G,$ allowed to have multiple edges but without
loops. We denote the vertex and edge set of $G$ by $V(G)$ and $E(G),$ respectively.
Given $u, v\in V(G),$ we set $a_{uv}$ to be equal to the number of edges between
vertices $u$ and $v.$ The matrix $A=A(G)=\{a_{uv}\}_{u, v\in V(G)}$ is called
\textit{the adjacency matrix} of the graph $G.$ The degree $d(v)$ of a vertex
$v \in V(G)$ is defined by $d(v)=\sum_{u\in V(G)}a_{uv}.$ Let $D=D(G)$ be the
diagonal matrix indexed by the elements of $V(G)$ with $d_{vv} = d(v).$ The matrix
$L=L(G)=D(G)-A(G)$ is called \textit{the Laplacian matrix}, or simply \textit{Laplacian},
of the graph $G.$

Recall \cite{Lor} the following useful relation between the structure of the
Laplacian matrix and the Jacobian of a graph $G.$ Consider the Laplacian $L(G)$
as a homomorphism $\mathbb{Z}^{|V|}\to\mathbb{Z}^{|V|},$ where $|V|=|V(G)|$ is
the number of vertices in $G.$ The cokernel
$\textrm{coker}\,(L(G))=\mathbb{Z}^{|V|}/\textrm{im}\,(L(G))$ is an abelian group. Let
$$\textrm{coker}\,(L(G))\cong\mathbb{Z}_{d_{1}}\oplus\mathbb{Z}_{d_{2}}\oplus\cdots\oplus\mathbb{Z}_{d_{|V|}}$$
be its Smith normal form satisfying the conditions $d_i\big{|}d_{i+1},\,(1\le i\le|V|-1).$
If the graph is connected, then the groups ${\mathbb Z}_{d_{1}},{\mathbb Z}_{d_{2}},\ldots,{\mathbb Z}_{d_{|V|-1}}$
are finite, and $\mathbb{Z}_{d_{|V|}}=\mathbb{Z}.$ In this paper, we define Jacobian
group $Jac(G)$ as the torsion subgroup of $\textrm{coker}\,(L(G)).$ In other words,
$$Jac(G)\cong\mathbb{Z}_{d_{1}}\oplus\mathbb{Z}_{d_{2}}\oplus\cdots\oplus\mathbb{Z}_{d_{|V|-1}}.$$

Let $M$ be an integer $n\times n$ matrix, then we can interpret $M$ as $\mathbb{Z}$-linear operator  from $\mathbb{Z}^n$ to $\mathbb{Z}^n.$ In this interpretation $M$ has a kernel
$\textrm{ker}M,$ an image $\textrm{im}\, M,$ and a cokernel
$\textrm{coker} M = \mathbb{Z}^n/\textrm{im}\, M.$ We emphasize that $\textrm{coker}\,M$ of the matrix
$M$ is completely determined by its Smith normal form. In particular, if  matrices $M$ and $M^{\prime}$  are elementary equivalent then $\textrm{coker}\,M\cong\textrm{coker}\,M^{\prime}.$

In what follows, by $I_n$ we denote the identity matrix of order $n.$

We call an $n\times n$ matrix {\it circulant,} and denote it by $circ(a_0, a_1,\ldots,a_{n-1})$ if it is of the form
$$circ(a_0, a_1,\ldots, a_{n-1})=
\left(\begin{array}{ccccc}
a_0 & a_1 & a_2 & \ldots & a_{n-1} \\
a_{n-1} & a_0 & a_1 & \ldots & a_{n-2} \\
  & \vdots &   & \ddots & \vdots \\
a_1 & a_2 & a_3 & \ldots & a_0\\
\end{array}\right).$$

Recall \cite{PJDav} that the eigenvalues of matrix $C=circ(a_0,a_1,\ldots,a_{n-1})$
are given by the following simple formulas $\lambda_j=p(\varepsilon^j_n),\,j=0,1,\ldots,n-1$
where $p(x)=a_0+a_1 x+\ldots+a_{n-1}x^{n-1}$ and $\varepsilon_n$ is an order $n$ primitive
root of the unity. Moreover, the circulant matrix $C=p(T),$ where $T=circ(0,1,0,\ldots,0)$
is the matrix representation of the shift operator
$T:(x_0,x_1,\ldots,x_{n-2},x_{n-1})\rightarrow(x_1, x_2,\ldots,x_{n-1},x_0).$
\smallskip

\section{The structure of Jacobian group for  the graph $Y(n;  k,l,m)$}

To investigate the   Laplacian matrix of graph $Y(n;  k,l,m)$, we denote by $T=circ(0,1,\ldots,0)$  the $n \times n$ shift operator.  Then the Laplacian $L=L(Y(n;  k,l,m))$  can be represented in the form
$$ L=\left(\begin{array}{cccc}
3I_n   & -I_n & -I_n & -I_n\\
-I_n & A & 0 & 0 \\
-I_n & 0 &B& 0 \\
-I_n & 0 & 0 & C\\
\end{array}\right),$$
where $A=3I_{n}-T^{-k}-T^k,B=3I_{n}-T^{-l}-T^l,\text{ and }C=3I_{n}-T^{-m}-T^m.$

Recall that $Jac(Y(n;k,l,m))$ is given by the torsion subgroup  of the $coker(L),$ where $L$ is considered as $\mathbb{Z}$-linear operator from $\mathbb{Z}^{4n}$ to itself. To find the structure of $coker(L)$ consider the following $n$-tuples of variables $x=(x_1,\ldots,x_n),\,y=(y_1,\ldots,y_n),\,s=(s_1,\ldots,s_n)$ and $t=(t_1,\ldots,t_n).$  As an abelian group,  $coker(L)$ has the following presentation
  $$coker(L)=\langle x,y,s,t|L(x,y,s,t)^t=0\rangle.$$  Hence,

$$coker(L)=\langle x,y,s,t|3x-y-s-t=0,\,-x+A y=0,\,-x+B s=0,\,-x+C t=0\rangle.$$
Since  $x=B s,$ we can eliminate the set of variables $x=(x_1,\ldots,x_n)$ from the above  presentation to get
$$coker(L)=\langle y,s,t|3B s-s-t-y=0,\,-B s+A y=0,\,-B s+C t=0\rangle.$$
Eliminating  $y=3B s -s-
t$ we also obtain
\begin{equation}\label{cok}coker(L)=\langle s,t|(-A-B+3A B)s-A t=0,\,-B s+C t=0\rangle .\end{equation}

This leads to the main result of this section.

\bigskip
\begin{theorem}\label{YgraphJac} The group $Jac(Y(n;k,l,m))$ is isomorphic to the torsion subgroup  of the $coker(M),$ where $M$ is the $2n\times2n$ block matrix given by
$$ M=\left(\begin{array}{cc}
-A-B+3A B  & -A \\
-B & C   \\
\end{array}\right),$$
where $A=3I_{n}-T^{-k}-T^k,\,B=3I_{n}-T^{-l}-T^l,\text{ and }C=3I_{n}-T^{-m}-T^m.$
\end{theorem}

\section{Explicit formulas for Jacobian of the graph $Y(n; 1,1,1)$}

We apply the results obtained in previous section to find explicit formulas for Jacobian of the graph $Y(n; 1,1,1).$ The main result of this section is  the following theorem.

\bigskip
\begin{theorem}\label{jacobian} Jacobian group $J_n=Jac(Y(n; 1,1,1)),\,n\ge 4$ of the graph $Y(n; 1,1,1)$  has the following structure
\begin{enumerate}
\item[ $1^0$ ]  $J_n\cong \mathbb{Z}_{3}^{n-4}\oplus\mathbb{Z}_{3n}\oplus\mathbb{Z}_{L_n}^2\oplus\mathbb{Z}_{3L_n}^2,$ if   $n$ is odd,
\item[ $2^0$ ]   $J_n\cong \mathbb{Z}_{3}^{n-4}\oplus\mathbb{Z}_{3n}\oplus\mathbb{Z}_{F_n}\oplus\mathbb{Z}_{3F_n}\oplus\mathbb{Z}_{5F_n}\oplus\mathbb{Z}_{15F_n},$  $n$ is even,
\end{enumerate}where $L_n$ and $F_n$  are the Lucas  and the Fibonacci numbers respectively.
\end{theorem}

The proof of the theorem is based on some auxiliary results. We start with the following result which is a direct consequence of Theorem~\ref{YgraphJac}.

\bigskip
\begin{proposition}\label{Y111Jac} The group $Jac(Y(n;1,1,1))$  is isomorphic to the torsion subgroup  of the abelian group
$coker(3(A^2-A))\oplus coker(A),$
where $A=3I_n-T -T^{-1}.$
\end{proposition}

\textbf{Proof} Since  $A=B=C=3I_n-T-T^{-1},$  by (\ref{cok})  we have

\begin{equation}\label{cok1}coker(L)=\langle s,t|(3A^2-2A)s-A t=0,\,-A s+A t=0\rangle .\end{equation}
Substituting   $t=s+u$ in (\ref{cok1}) we obtain
\begin{eqnarray}\label{cok2}
\nonumber coker(L)&=&\langle s,u|(3A^2-3A)s-A u=0,\,A u =0\rangle=\langle s,u|(3(A^2-A)s=0,\,A u=0\rangle \\
\nonumber &=&\langle s|(3(A^2-A)s=0\rangle\oplus\langle u|\,A u=0\rangle=coker(3(A^2-A))\oplus coker(A).
\end{eqnarray}
\

Also, we need the following version of Theorem 2.1(iv) from \cite{CorValen}.

\bigskip
\begin{proposition}\label{a_plus_b} Let $T=circ(0,1,0,\ldots,0)$ be the $n$ by $n$ shift operator and $a$ and $b$
be non-zero integers. We set $C_{n}(a,b)=aI_{n}+b T+b T^{-1}.$ Then $coker\,C_{n}(a,b)\cong coker\,C,$ where

$$C=\begin{cases}
 b^{\frac{n}{2}-1}u(n)\left(\begin{array}{cc}
a & 2b\\
2b & a\\
\end{array}\right) & \textrm{if } $n$ \textrm{  is even }\\
 b^{\frac{n-1}{2}}t(n)
\left(\begin{array}{cc}
1  & 0 \\
0 & a+2b \\
\end{array}\right) & \textrm{if } $n$ \textrm{  is odd,}
\end{cases}
$$
\end{proposition}
   $$ u(n)=\frac{\sin(\frac{n}{2}\arccos(x))}{\sin(\arccos(x))},\,\,
t(n)=\frac{\cos(\frac{n}{2}\arccos(x))}{\cos(\frac{1}{2}\arccos(x))}\text{ and }x=\frac{a}{2b}.$$

\textbf{Proof} Consider a two-variable Chebyshev polynomial $f_n(x, y)$ defined by the following  recursion
$$f_n(x, y)=x f_{n-1}(x, y)-y^2f_{n-2}(x, y),\,f_{-1}(x, y)=1,\,f_0(x, y)=0.$$ One can check that $f_n(x, y)=y^n U_n(\frac{x}{2y}),$ where $U_n(z)$ is the Chebyshev polynomial of the second kind. By (\cite{CorValen}, Theorem 2.1(iv)) we have $coker\,C_{n}(a,b)\cong coker\,C,$ where

$$C=\begin{cases}
f_{\frac{n}{2}-1}(a,b)\left(\begin{array}{cc}
a & 2b\\
2b & a\\
\end{array}\right) & \textrm{if } n \textrm{  is even}\\
 (f_{\frac{n-1}{2}}(a,b)-b f_{\frac{n-3}{2}}(a,b))
\left(\begin{array}{cc}
1  & 0 \\
0 & a+2b \\
\end{array}\right) & \textrm{if } n  \textrm{  is odd,}
\end{cases}
$$ To finish the proof we have to show that  $f_{\frac{n}{2}-1}(a,b)= b^{\frac{n}{2}-1}u(n)$  and
$ f_{\frac{n-1}{2}}(a,b)-b f_{\frac{n-3}{2}}(a,b)=b^{\frac{n}{2}-1}t(n),$
where $u(n)$ and $t(n)$ are the same as in the statement of Proposition.  The first equality easily    follows from the definition of   Chebyshev  polynomial, since  $$U_{\frac{n}{2}-1}(x)=\frac{\sin(\frac{n}{2}\arccos(x))}{\sin(\arccos(x))}.$$
 To prove the second we note  that  $f_{\frac{n-1}{2}}(a,b)=b^{\frac{n-1}{2}}U_{\frac{n-1}{2}}(x)=b^{\frac{n-1}{2}}\frac{\sin(\frac{n+1}{2}\arccos(x))}{\sin(\arccos(x))}.$ Hence,
$$f_{\frac{n-1}{2}}(a,b)-b f_{\frac{n-3}{2}}(a,b)=b^{\frac{n-1}{2}}\left(\frac{\sin(\frac{n+1}{2}\arccos(x))}{\sin(\arccos(x))}-\frac{\sin(\frac{n-1}{2}\arccos(x))}{\sin(\arccos(x))}\right)$$

$$=2 b^{\frac{n-1}{2}}\left(\frac{\cos(\frac{n}{2}\arccos(x))\sin(\frac{1}{2}\arccos(x))}{\sin(\arccos(x))}\right)=b^{\frac{n-1}{2}}\frac{\cos(\frac{n}{2}\arccos(x))}{\cos(\frac{1}{2}\arccos(x))}.$$

The following two lemmas are immediate consequences of Proposition~\ref{a_plus_b}.
\bigskip
\begin{lemma}\label{threeone}
Let $T=circ(0,1,0,\ldots,0)$  be the $n$ by $n$  shift operator. We  set  $A(n)=3I_n- T^{-1}-T.$  Then cokernel of  $A(n)$  is given by the following formulas
\begin{enumerate}
\item[(i)] $coker\,A(n)\cong\mathbb{Z}_{F_n}\oplus\mathbb{Z}_{5F_n}$  if  $n$ is even,
\item[(ii)] $coker\,A(n)\cong\mathbb{Z}_{L_n}\oplus\mathbb{Z}_{L_n}$  if  $n$ is odd,
\end{enumerate}
 where $F_n$  and $L_n$  are the Fibonacci and the  Lucas numbers respectively.
\end{lemma}

\textbf{Proof} Since $A(n)=C_n(3,-1),$ we have $coker\,A(n)\cong coker(C),$ where $C$ is the same as in Proposition~\ref{a_plus_b}, with $a=3$  and $b=-1.$
Then  for $n$  even we get $$b^{\frac{n}{2}-1}u(n)=(-1)^{\frac{n}{2}-1}\frac{\sin(\frac{n}{2}\arccos(-\frac{3}{2}))}{\sin(\arccos(-\frac{3}{2}))}=\frac{1}{\sqrt{5}}\left((\frac{\sqrt{5}+1}{2})^n-(\frac{\sqrt{5}-1}{2})^n\right)=F_n,$$
and for $n$  odd

$$b^{\frac{n-1}{2}}t(n)=(-1)^{\frac{n-1}{2}}\frac{\cos(\frac{n}{2}\arccos(-\frac{3}{2}))}{\cos(\frac{1}{2}\arccos(-\frac{3}{2}))}= (\frac{\sqrt{5}+1}{2})^n-(\frac{\sqrt{5}-1}{2})^n=L_n.$$  Hence,  $C= \left(\begin{array}{cc}
3 F_n & -2F_n \\
-2F_n & 3F_n \\
\end{array}\right) $  if $n$  is even   and  $C= \left(\begin{array}{cc}
L_n  & 0 \\
0 & L_n \\
\end{array}\right) $  if $n$  is odd.  We note that matrix $\left(\begin{array}{cc}
3  & -2 \\
-2 & 3 \\
\end{array}\right)$ is elementary equivalent to matrix  $\left(\begin{array}{cc}
1  & 0 \\
0 & 5 \\
\end{array}\right).$
Therefore, $coker(C)\cong\mathbb{Z}_{F_n}\oplus\mathbb{Z}_{5F_n}$ if $n$ is even and
$coker(C)\cong\mathbb{Z}_{L_n}\oplus\mathbb{Z}_{L_n}$ if $n$ is odd. The lemma is proved.

By similar arguments we get the following result.
\bigskip
\begin{lemma}\label{twoone}
Let $T=circ(0,1,0,\ldots,0)$ and $B(n)= 2I_n-T^{-1}-T .$ Then for $n\ge2$ we have
$$coker\,B(n)\cong \mathbb{Z}_{n}\oplus\mathbb{Z}.$$
\end{lemma}

\bigskip

We also need the following elementary lemma.

\begin{lemma}\label{three} Let  $D={diag}(d_1,d_2,\ldots,d_n)$ be an integer $n\times n$
diagonal matrix  and $k$  be a positive integer.  Then
$$coker(k D)\cong \mathbb{Z}_{k d_1}\oplus\mathbb{Z}_{k d_2}
\ldots\oplus\mathbb{Z}_{k d_n}.$$
\end{lemma}

\bigskip

Now, we are able to complete the proof of Theorem~\ref{jacobian}. Indeed, by Proposition~\ref{Y111Jac} $coker(L)=coker(3(A^2-A))\oplus coker(A).$ The structure of $coker(A)$ is described in Lemma~\ref{threeone}.
Also, Lemma~\ref{twoone} gives us $coker(A-I)\cong\mathbb{Z}_{n}\oplus\mathbb{Z},$  where $I=I_n.$

We note that the matrices
$\left(\begin{array}{cc}
A  & 0  \\
I & A-I \\
\end{array}\right)$
and
$\left(\begin{array}{cc}
A  & 0  \\
0 & A-I \\
\end{array}\right)$ are elementary equivalent.

Hence, $coker(A(A-I))\cong coker(A)\oplus coker(A-I).$ By Lemmas \ref{threeone}  and \ref{twoone},
$coker(A(A-I))\cong\mathbb{Z}_{F_n}\oplus\mathbb{Z}_{5F_n}\oplus\mathbb{Z}_{n}\oplus\mathbb{Z}$  if $n$ is even and
$coker(A(A-I))\cong\mathbb{Z}_{L_n}\oplus\mathbb{Z}_{L_n}\oplus\mathbb{Z}_{n}\oplus\mathbb{Z}$  if $n$ is odd. This means that matrix $A(A-I)$   is elementary equivalent to the matrix 
$$ {diag}(\underbrace{1,\ldots,1}_{(n-4)\, times},F_n, 5F_n, n, 0)$$  
in the first case and 
to  the matrix
$$ {diag}(\underbrace{1,\ldots,1}_{ (n-4)\, times },L_n, L_n, n,0)$$  
in the latter. Then,
by making use of Lemma~\ref{three}, we  get
$coker(3A(A-I))\cong\mathbb{Z}_{3}^{n-4}\oplus\mathbb{Z}_{3F_n}\oplus\mathbb{Z}_{15F_n}\oplus\mathbb{Z}_{3n}\oplus\mathbb{Z}$  if $n$ is even and
$coker(3A(A-I))\cong\mathbb{Z}_{3}^{n-4}\oplus\mathbb{Z}_{3L_n}\oplus\mathbb{Z}_{3L_n}\oplus\mathbb{Z}_{3n}\oplus\mathbb{Z}$  if $n$ is odd. Once more, by Lemma~\ref{threeone} we have $coker(A)\cong\mathbb{Z}_{F_n}\oplus\mathbb{Z}_{5F_n}$
 if $n$ is even and   $coker(A)\cong\mathbb{Z}_{L_n}\oplus\mathbb{Z}_{L_n}$
 if $n$ is odd.
 
Combine these results we finish the proof of the theorem.

\section{Counting the number  of spanning trees for the graph $Y(n;k,l,m)$}

\begin{theorem}\label{theorem1} The number of spanning trees $\tau_{k,l,m}(n)$ in the   graph $Y(n;k,l,m)$
is given by the formula $$\tau_{k,l,m}(n)=\frac{n\,3^{n}}{k^2+l^2+m^2 }
\prod_{p=1}^{s}|2T_n(w_p)-2|,$$ hereby $
s=k+l+m-1, $ and $\, w_p=1,2,\ldots,s$ are roots of the order $s$ algebraic equation
$Q(w)=0,$ where
$$Q(w)=\frac{(3-2T_k(w))(3-2T_l(w))(3-2T_m(w))}{w-1}\left(3-\frac{1}{3-2T_k(w)}-\frac{1}{3-2T_l(w)}-\frac{1}{3-2T_m(w)}\right)$$
and $T_k(w)$ is the Chebyshev
polynomial of the first kind.
\end{theorem}

\medskip\noindent{Proof:} \ \
        By the celebrated Kirchhoff theorem, the number of spanning trees $\tau_{k,\,l,\,m}(n)$ is equal to the
product of nonzero eigenvalues of the Laplacian of a graph $Y(n;k,l,m)$ divided by the number of its vertices $4n.$
To investigate the spectrum of Laplacian matrix, we denote by $T=circ(0,1,\ldots,0)$
  the $n \times n$ shift operator.  The the Laplacian $L=L(Y(n;k,l,m))$  is given by the matrix
$$ L=\left(\begin{array}{cccc}
3I_n   & -I_n & -I_n & -I_n\\
-I_n & 3I_n-T^{-k}-T^{k}& 0 & 0 \\
-I_n & 0 &3I_n-T^{-l}-T^{l}& 0 \\
-I_n & 0 & 0 & 3I_n-T^{-m}-T^{m}\\
\end{array}\right).$$

The eigenvalues of circulant matrix $T$ are $\varepsilon_n^j,\,j=0,1,\ldots,n-1,$ where
$\varepsilon_n=e^\frac{2\pi i}{n}.$ Since all of them are
distinct, the matrix $T$ is conjugate to the diagonal matrix
$\mathbb{T}=diag(1,\varepsilon_n,\ldots,\varepsilon_n^{n-1})$ with
diagonal entries
$1,\varepsilon_n,\ldots,\varepsilon_n^{n-1}$. To find spectrum of
$L,$ without loss of generality, one can assume that $T=\mathbb{T}.$
Then the $n \times n$ blocks of $L$ are diagonal matrices. This
essentially simplifies the problem of finding eigenvalues of $L.$
Indeed, let $\lambda$ be an eigenvalue of $L$ and
$(x,y,z,t)=(x_1,\ldots,x_n,y_1,\ldots,y_n,z_1,\ldots,z_n,t_1,\ldots,t_n)$ be the respective
eigenvector. Then we have the following system of equations

$$\left\{\begin{array}{ccccc}
(3-\lambda)x &-y &-z&-t &=0\\
-x&+(A -\lambda)y &{}&{} &= 0 \\
-x&{}&+(B -\lambda)z  &{}&= 0 \\
-x&{} &{}&+(C -\lambda)t &= 0 \\
\end{array}\right.,$$
where $A=3I_{n}-T^{-k}-T^k,\,B=3I_{n}-T^{-l}-T^l,\text{ and }C=3I_{n}-T^{-m}-T^m.$


Consequently eliminating $y, z$  and $t$  from the above equation, we obtain
$$((3-\lambda)(A-\lambda)(B-\lambda)(C-\lambda) -(A-\lambda)(B-\lambda)-(B-\lambda)(C-\lambda)-(A-\lambda)(C-\lambda) )x=0.$$
Recall that the matrices under consideration are diagonal and the $(j+1,j+1)$-th entry of $T$ is equal to $\varepsilon_n^{j}.$  Consider the following two variable Laurent polynomial
$$P(z,\lambda)=(3-\lambda)(a-\lambda)(b-\lambda)(c-\lambda) -(a-\lambda)(b-\lambda)-(b-\lambda)(c-\lambda)-(a-\lambda)(c-\lambda),$$ where
$a=3-z^k-z^{-k},\, b=3-z^l-z^{-l}$  and $c=3-z^m-z^{-m}.$

Then,  for any $j=0,\ldots, n-1,$   matrix $L$ has four  eigenvalues,
say $\lambda_{1,j},\lambda_{2,j},\lambda_{3,j},$ and $\lambda_{4,j},$ satisfying the order three algebraic equation
$ P(\varepsilon_n^{j},\lambda)=0.$   In particular, if $j=0$ for $\lambda_{1,0}, \lambda_{2,0},\lambda_{3,0}, \lambda_{4,0}$ we have $\lambda(\lambda-1)^2(\lambda-4)=0.$
Thus, $\lambda_{1,0}=0,\lambda_{2,0}=1,\lambda_{3,0}=1,$ and $\lambda_{4,0}=4.$

Let ${\bf e}_{j} =(0,\ldots,\underbrace{1}_{j-th},\ldots, 0),\,j=1,\ldots, n.$
For given  $\lambda=\lambda_{i,j}$ the corresponding eigenvector  is  $(x,y,z,t),$ where $x={\bf e}_{j+1}, y= (B-\lambda)(C-\lambda){\bf e}_{j+1},z= (A-\lambda)(C-\lambda){\bf e}_{j+1},$ and $t=(A-\lambda)(B-\lambda){\bf e}_{j+1}.$

By   Vieta's theorem, the product  $\pi_{j}=\lambda_{1,j}\lambda_{2,j}\lambda_{3,j}\lambda_{4,j}$ is given by the formula $\pi_{j}=P(\varepsilon_n^j),$  where  $P(z)=P(z,0).$

Now we have $$\tau_{k,\,l,\,m}(n)=\frac{1}{4n}\lambda_{2,0}\lambda_{3,0}\lambda_{4,0}\prod\limits_{j=1}^{n-1}\lambda_{1,j}\lambda_{2,j}\lambda_{3,j}\lambda_{4,j}=
\frac{1}{n}\prod\limits_{j=1}^{n-1}\pi_j=\frac{1}{n}\prod\limits_{j=1}^{n-1}P(\varepsilon_n^j).$$

To continue the calculation of $\tau_{k,\,l,\,m}(n)$ we need the following two lemmas.

\begin{lemma}\label{lemma2} The following identity holds
$P(z)=(w-1)Q(w),$  where  $Q(w)$ is an order $k+l+m-1$ integer polynomial given by the formula
 $$Q(w)=\frac{(3-2T_k(w))(3-2T_l(w))(3-2T_m(w))}{w-1}\left(3-\frac{1}{3-2T_k(w)}-\frac{1}{3-2T_l(w)}-\frac{1}{3-2T_m(w)}\right),$$
  $T_k(w)$ is the Chebyshev polynomial of the first kind and
$w=\frac{1}{2}(z+z^{-1}).$
\end{lemma}

\textbf{Proof}  Recall that $P(z)=3a b c- a b - bc -ac,$ where $a=3-z^k-z^{-k},\, b=3-z^l-z^{-l}$  and $c=3-z^m-z^{-m}.$ Let us substitute $z=e^{i\varphi}.$ It is easy to
see that $w=\frac{1}{2}(z+z^{-1})=\cos\varphi,$ so we have
$T_k(w)=\cos(k\arccos w)=\cos(k\varphi)=\frac{1}{2}(z^k+z^{-k})$. Then the statement of the
lemma follows from elementary calculations.

Note that the leading term of Laurent polynomial $P(z)$ is $-3 z^{k+l+m}.$
By Lemma \ref{lemma2}, $P(z)=(w-1)Q(w),$ where $w=\frac{1}{2}(z+z^{-1})$ and $Q(w)$ is a polynomial of degree $s=k+l+m-1$. We have $w-1=\frac{(z-1)^2}{2z}$.  Since $Q(1)=\lim_{w\to1}Q(w)=-2(k^2+l^2+m^2)\neq 0,$   polynomial $P(z)$ has the root $z=1$ with
multiplicity two. Hence, all the roots of $P(z)$ are
$1,1,z_1,1/z_1,\ldots,z_s,1/z_s,$ where for all $p=1,\ldots,s,
z_p\neq 1.$ We get
$$P(z)=\frac{-3(z-1)^2}{z^{k+l+m}}\prod\limits_{p=1}^{s}(z-z_p)(z-z_p^{-1}).$$
Also,  $w_p=\frac{1}{2}(z_p+z_p^{-1}),\,p=1,\ldots,s$ are all the roots of the  equation
$Q(w)=0.$ We set $H(z)=\prod\limits_{p=1}^{s}(z-z_p)(z-z_p^{-1}).$
Then $P(z)=\frac{-3(z-1)^2}{z^{k+l+m}}H(z).$

\medskip
The following lemma has been proved in (\cite{KMM}, Lemma 5.3).
\begin{lemma}\label{lemma5}  Let $H(z)=\prod\limits_{p=1}^{s}(z-z_p)(z-z_p^{-1})$ and $H(1)\neq0$. Then
$$\prod\limits_{j=1}^{n-1}H(\varepsilon_n^j)=\prod\limits_{p=1}^{s}\frac{T_n(w_p)-1}{w_p-1},$$  where
$w_p=\frac12(z_p+z_p^{-1}),\,p=1,\ldots,k$  and $T_n(x)$  is the
Chebyshev polynomial of the first kind.
\end{lemma}

Note that
$\prod\limits_{j=1}^{n-1}(\varepsilon_n^j-1)^2=\lim\limits_{z\to1}\prod\limits_{j=1}^{n-1}
(z-\varepsilon_n^j)^2= \lim\limits_{z\to1}(\frac{z^n-1}{z-1})^2=n^2$ and
$\prod\limits_{j=1}^{n-1}\varepsilon_n^{j} = (-1)^{n-1}$.  As a
result, taking into account Lemma~\ref{lemma5}, we obtain
\begin{eqnarray*}
\tau_{k,l,m}(n)&=&\frac{1}{n}\prod\limits_{j=1}^{n-1}P(\varepsilon_n^j)=\frac{1}{n}\prod\limits_{j=1}^{n-1}
\frac{(-3)(\varepsilon_n^j-1)^2}{(\varepsilon_n^{j})^{s+1}}H(\varepsilon_n^j)\\
&=&\frac{(-1)^{s(n-1)}3^{n-1}n^2}{n}
\prod\limits_{j=1}^{n-1}H(\varepsilon_n^j)=
(-1)^{s(n-1)}3^{n-1}n\prod\limits_{p=1}^{s}\frac{T_n(w_p)-1}{w_p-1}.
\end{eqnarray*}

We note that $w_p,\, p=1,\ldots,s$ are  all the roots of the polynomial  $Q(w)$ is an order $s=k+l+m-1$   whose  leading term is $-3\cdot2^{k+l+m}.$  We also have $Q(1)=-2(k^2+l^2+m^2).$ Hence

$$\prod\limits_{p=1}^{s}(w_p-1)=\frac{(-1)^s Q(1)}{-3\cdot2^{k+l+m}}=\frac{k^2+l^2+m^2}{3\cdot(-2)^{s}}.$$  Finally,

$$\tau_{k,l,m}(n)=(-1)^{n s}\frac{n\,3^{n}2^s}{k^2+l^2+m^2 }
\prod_{p=1}^{s}(T_n(w_p)-1)=\frac{n\,3^{n}}{k^2+l^2+m^2 }
\prod_{p=1}^{s}|2T_n(w_p)-2|.$$


\section{Asymptotic formulas for the number of spanning trees $\tau_{k,l,m}(n)$ }
Consider Laurent polynomial $P(z)=3abc-ab-ac-bc,$ where $a=3-z^k-z^{-k},\,b=3-z^l-z^{-l}$ and $c=3-z^m-z^{-m}$ from the above section. We get the following result.
\begin{theorem}\label{asymptotic}
The asymptotic behaviour  for number of spanning trees $\tau_{k,l,m}(n)$ for the graph $Y(n; k,l,m)$ with $\gcd(k,l,m)=1$ is given by the formula
$$\tau_{k,l,m}(n)\sim \frac{n}{k^2+l^2+m^2}A_{k,l,m}^n, \, n\to \infty,$$ where
$A_{k,l,m}=|a_{0}|\prod\limits_{P(z)=0,\,|z|>1}|z|,$ where the product is taken over roots of  the Laurent polynomial $P(z)$ and $a_{0}$ is the leading coefficient of $P(z).$
\end{theorem}

We note that  number $A_{k,l,m}$ coincides with the Mahler measure of polynomial $P(z).$ Hence

$$A_{k,l,m}=\exp\left({\int\limits_{0}^{1}\log|P(e^{2 \pi i t})|\textrm{d}t}\right).$$

\textbf{Proof}
By Theorem \ref{theorem1} we have
$\tau_{k,l,m}(n)=\frac{n\,3^{n}}{k^2+l^2+m^2}\prod_{p=1}^{s}|2T_n(w_p)-2|,$ where
$w_{j},\,j=1,2,\ldots,s$ are all the roots of the polynomial $Q(w).$
Recall $P(z)=(w-1)Q(w),$ where $w=\frac{1}{2}(z+\frac{1}{z}).$ All the roots of the polynomial $(w-1)Q(w)$ are $w_{j}=\frac{1}{2}(z_{j}+z_{j}^{-1}),$ where $z_{j}$ and $1/z_{j},\,j=1,\ldots,s+1$ are all the roots of the polynomial $P(z).$
 One can check that $P(z)$ has $1$ as a root of multiplicity two. Indeed, $P(1)=P^{\prime}(1)=0$ and $P^{\prime\prime}(1)=-2(k^2+l^2+m^2)\neq0.$ So, we can assume that
$z_{j}\neq1,\,j=1,2,\ldots,s.$ Since $\gcd(k,l,m)=1,$ taking into account the arguments from Lemma 5.2 in \cite{Ilya} we get $|z_{j}|\neq1,\,j=1,2,\ldots,s.$

We have $T_{n}(w_{j})=\frac{1}{2}(z_{j}^{n}+z_{j}^{-n}).$ Replacing
$z_{j}$ by $1/z_{j},$ if necessary, we always can assume that
$|z_j|>1$ for all $j=1,2,\ldots,s.$ Then
$T_{n}(w_{j})\sim\frac{1}{2}z_{j}^{n}$ and
$|2T_{n}(w_{j})-2|\sim|z_{j}|^{n}$ as $n\to\infty.$ Hence
$$3^{n}\prod_{j=1}^{s}|2T_{n}(w_{j})-2|\sim 3^{n}\prod_{j=1}^{s}|z_{j}|^{n}
=3^{n}\prod\limits_{P(z)=0,\,|z|>1}|z|^{n}=A_{k,l,m}^{n}.$$
The leading coefficient $a_0$ of the polynomial $P(z)$ is $3.$ So
$A_{k,l,m}=3\prod\limits_{P(z)=0,\,|z|>1}|z|$ is the Mahler measure of
the polynomial $P(z).$ By (\cite{EverWard}, p.~67), we have
$A_{k,l,m}=\exp\left(\int\limits_{0}^{1}\log|P(e^{2 \pi i t})|\textrm{d}t\right).$
The theorem is proved.

\section{Examples}

\subsection{Graph $Y(n,1,1,1)$}
Note that number $\tau_{1,1,1}(n)$ of spanning trees in
the   graph $Y(n,1,1,1)$ coinsides with the size of its Jacobian  $J_n.$
As a consequence of theorem \ref{jacobian} we obtain
\begin{enumerate}
\item[ $1^0$ ]  $\tau_{1,1,1}(n)=3^{n-1}n \,L_n^4, $ if   $n$ is odd,
and
\item[ $2^0$ ] $\tau_{1,1,1}(n)=25\cdot3^{n-1}n\, F_n^4,$ if   $n$ is even,
\end{enumerate}where  $L_n$  and $F_n$  are the Lucas  and the Fibonacci numbers respectively.

\subsection{Graph $Y(n,1,1,2)$} The number $\tau(n)=\tau_{1,1,2}(n)$ of spanning trees in
the   graph $Y(n,1,1,2)$  is given by the formula
$$\tau(n)={4n\,3^{n-1}}\left|\left(T_n({3\over2})-1\right)\left(T_n({1 + \sqrt{193}\over 12})-1\right)\left(T_n({1 - \sqrt{193}\over 12})-1 \right)\right|.$$

By Theorem~\ref{asymptotic}, $\tau_{1,1,2}(n)\sim \frac n6 A_{1,1,2}^n, \, n\to \infty,$  where  $$A_{1,1,2}=\frac{1}{4}(3+\sqrt{5})(8 + \sqrt{3} + \sqrt{31 + 16 \sqrt{3}})=22.7697\ldots.$$

\subsection{Graph $Y(n,1,2,2)$} We have the following formulas for number $\tau(n)=\tau_{1,2,2}(n)$ of spanning trees in
the   graph $Y(n,1,2,2).$
$$\tau(n)=n\,3^{n-2}\left((\frac{1 + \sqrt{13}}{2})^{n} +  (\frac{1 + \sqrt{13}}{2})^{-n}-2\right)\left((\frac{1 - \sqrt{13}}{2})^{n} + (\frac{1-\sqrt{13}}{2})^{-n}-2\right)$$

$$\times\left((\frac{1+\sqrt{5}}{2})^{n} + (\frac{1+\sqrt{5}}{2})^{-n}-2\right)\left((\frac{-1-\sqrt{5}}{2})^{n} + (\frac{-1-\sqrt{5}}{2})^{-n}-2\right).$$

Also, by  Theorem~\ref{asymptotic}, $\tau_{1,2,2}(n)\sim \frac n9 A_{1,2,2}^n, \, n\to \infty,$  where  $A_{1,2,2}=\frac92 (3 + \sqrt{5})=23.5623\ldots .$

One can show that there exists  an integer sequence $a(n)$ such that
 $\tau(n) = na(n)^2.$

\section{Final Remarks}

Theorem~\ref{YgraphJac} is the first step to understand the structure of the Jacobian for $Y(n;k,l,m).$  Also, it gives a way for calculations of $\textrm{Jac}(Y(n;k,l,m)$ for small
values of $n, k, l$ and $m.$ Theorem~\ref{jacobian} completely describes Jacobian for $Y(n;1,1,1).$ Theorems \ref{theorem1} and \ref{asymptotic}  contain convenient formulas for counting the number of spanning trees and  their  asymptotics.

\section{Acknowledgment}
The first author was supported  by the Basic Science Research Program through the National Research Foundation of Korea (NRF) funded by the Ministry of Education (2018R1D1A1B05048450). The second and the third authors were supported by the Mathematical Center in Akademgorodok under agreement
No.\,075-15-2019-1613 with the Ministry of Science and Higher Education of the Russian Federation.

\vspace{1cm}
\noindent Young  Soo Kwon

\noindent Department of Mathematics,  

\noindent Yeungnam University,

\noindent Gyeongsan, Gyeongbuk 38541, Korea%

\noindent \text{ysookwon@ynu.ac.kr}%

\vspace{0.75cm}
\noindent Aleksandr Dmitrievich  Mednykh

\noindent Sobolev Institute of Mathematics,

\noindent pr-t Acad. Koptyug, 4,

\noindent 630090, Novosibirsk, Russia\vspace{5pt}

\noindent Novosibirsk State University,

\noindent Pirogova str., 1,

\noindent 630090, Novosibirsk, Russia%

\noindent \text{smedn@mail.ru}%

\vspace{0.75cm}
\noindent Ilya Aleksandrovich Mednykh

\noindent Sobolev Institute of Mathematics,

\noindent pr-t Acad. Koptyug, 4,

\noindent 630090, Novosibirsk, Russia\vspace{5pt}

\noindent Novosibirsk State University,

\noindent Pirogova str., 1,

\noindent 630090, Novosibirsk, Russia%

\noindent \text{ilyamednykh@mail.ru}%

\end{document}